\documentclass[final]{elsarticle}

\usepackage{lineno,hyperref}\usepackage{amsmath,amssymb}
\usepackage[color]{showkeys}
\definecolor{refkey}{rgb}{0,1,1}
\definecolor{labelkey}{rgb}{1,0,0}
\modulolinenumbers[100]
\usepackage{graphicx}

\makeatletter
\def\ps@pprintTitle{%
  \let\@oddhead\@empty
  \let\@evenhead\@empty
  \def\@oddfoot{\reset@font\hfil\thepage\hfil\today}
  \let\@evenfoot\@oddfoot
}
\makeatother

\newtheorem{thm}{Theorem}

\newdefinition{rmk}{Remark}
\newproof{pf}{Proof}
\newproof{pot}{Proof of Theorem \ref{thm2}}
\newcommand{\norm}[1]{\left\Vert#1\right\Vert}
\newcommand{\abs}[1]{\left\vert#1\right\vert}
\newcommand{\eq} [1] {\begin{equation}\label{#1}\quad}
\newcommand{\en} {\end{equation}}
\newcommand{\es}{\operatorname{ess\,sup}}
\newcommand{\dg}{\operatorname{diag}}
\renewcommand{\Re}{\operatorname{Re}}
\begin{document}

\begin{frontmatter}

\title{A distance formula related to a family of projections orthogonal to their symmetries}

\author{Ilya M. Spitkovsky}
\address{Division of Science,  New York  University Abu Dhabi (NYUAD)\\ Saadiyat Island,
P.O. Box 129188 Abu Dhabi, UAE}
\begin{abstract}
Let $u$ be a hermitian involution,  and $e$ an orthogonal projection, acting on the same Hilbert space $\mathcal H$. We establish the exact formula, in terms of $\norm{eue}$, for the distance from $e$ to the set of all orthogonal projections $q$ from the algebra generated by $e,u$, and such that $quq=0$.
\end{abstract}

\begin{keyword}
orthogonal projection \sep involution \sep $C^*$-algebra \sep $W^*$-algebra
\MSC[2010] 47A05 \sep  47A30
\end{keyword}

\end{frontmatter}

\linenumbers

\section{Introduction}

Let $\mathcal H$ be a  Hilbert space and let $\mathcal{B}(\mathcal{H})$ stand for the $C^*$-algebra of all bounded linear operators acting on $\mathcal{H}$. Given a hermitian involution $u\in\mathcal{B}(\mathcal{H})$,  denote by ${\mathcal Q}_u$ the set of all orthogonal projections $q\in\mathcal{B}(\mathcal{H})$ for which $quq=0$.

Theorem~1.2 of \cite{Walt17} can be stated as follows:
\begin{thm}\label{th:Wal} Let $e\in\mathcal{B}(\mathcal{H})$ be an orthogonal projection such that \[ \norm{eue}<\xi (\approx 0.455).\]  Then there exists $q\in{\mathcal Q}_u$ for which
\eq{eq} \norm{e-q}\leq\frac{1}{2}\norm{eue}+4\norm{eue}^2.\en
Further, $q$ is in the $C^*$-subalgebra of $\mathcal{B}(\mathcal{H})$ generated by $e,ueu^*$.
\end{thm}
Note that the distance between any two orthogonal projections does not exceed one. So, estimate \eqref{eq} is useful only when $\norm{eue}$ is smaller than the positive root of $8x^2+x-2$, that is, approximately $0.441$.

We will provide an explicit formula for the distance from $e$ to the intersection of ${\mathcal Q}_u$ with the $W^*$-algebra ${\mathcal W}(e,u)$ generated by $e,u$, as well as for the
element on which this distance is attained. No a priori restriction on $\norm{eue}$ is needed, and the respective $q$ indeed lies in the $C^*$-algebra ${\mathcal C}(e,ueu^*)$ generated by $e,ueu^*$ whenever $\norm{eue}<1$.

\begin{thm}\label{th:main}Let $e,u\in\mathcal{B}(\mathcal{H})$ be, respectively, an orthogonal projection and a hermitian involution. Denote by ${\mathcal H}_\pm$ the eigenspace of $u$ corresponding to
its eigenvalue $\pm 1$. Then the distance $d$ from $e$ to ${\mathcal Q}_u\cap{\mathcal W}(e,u)$ is one
if the range of $e$ has a non-trivial intersection with ${\mathcal H}_+$ or ${\mathcal H}_-$, and is given by the formula
\eq{d} d= \sqrt{\frac{1}{2}\left(1-\sqrt{1-\norm{eue}^2}\right)} \en otherwise. \end{thm}
For small values of $\norm{eue}$, it is instructive to compare \eqref{eq} with the Taylor expansion of \eqref{d}:
\[ d=\frac{1}{2}\norm{eue}+\frac{1}{16}\norm{eue}^3+\cdots \]
\begin{figure}[h]
\center{\includegraphics[scale=0.6]{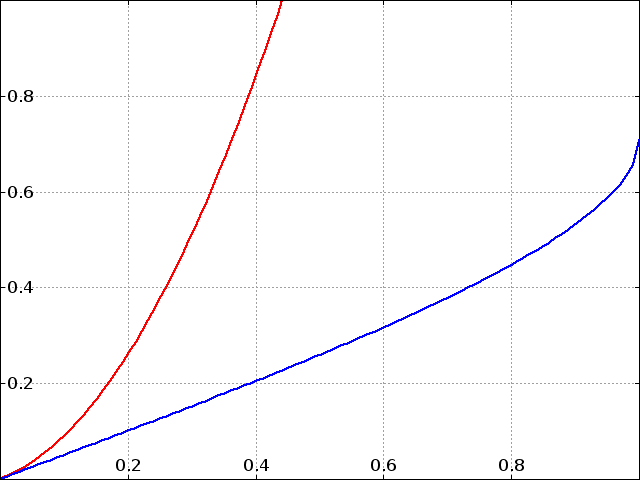}}
\caption{Estimate \eqref{eq} versus formula \eqref{d} as functions of $\norm{eue}$}
\end{figure}

\section{Proof of the main result} \label{s:pr}
Using the canonical representation \cite{Hal69} (see also \cite{Spit94} or a more recent survey \cite{BSpit10}) of the  pair $e, (u+I)/2$ of orthogonal projections, we can find an orthogonal
decomposition of $\mathcal H$ into six summands,
\eq{H} {\mathcal H}={\mathcal M}_{00}\oplus{\mathcal M}_{01}\oplus{\mathcal M}_{10}\oplus{\mathcal M}_{11}\oplus\left({\mathcal M}\oplus{\mathcal M}\right), \en
with respect to which
\eq{tp} \begin{split}  u = &  I  \oplus  I  \oplus  (-I)  \oplus   (-I)  \oplus  \dg [I, -I], \\
 e = & I  \oplus  0 \oplus  0  \oplus I  \oplus  \begin{bmatrix}  H & \sqrt{H(I-H)} \\ \sqrt{H(I-H)} & I-H \end{bmatrix}. \end{split}  \en
(Here and in what follows we use the notation $\dg[X_1,\ldots,X_k]$ for block diagonal matrices with $X_1,\ldots X_k$ as their diagonal blocks.)
Note that in \eqref{H} the subspaces ${\mathcal M}_{00}$ and ${\mathcal M}_{11}$ (resp, ${\mathcal M}_{01}$ and ${\mathcal M}_{10}$) are the intersections of the range (resp, the kernel) of $e$ with ${\mathcal H}_+$
and ${\mathcal H}_-$. The (hermitian) operator $H$ is the compression of $e$ onto ${\mathcal M}:={\mathcal H}_+\ominus({\mathcal M}_{00}\oplus{\mathcal M}_{01})$. By construction, $H$
has its spectrum $\Delta$ lying in $[0,1]$ and $0,1$ are not its eigenvalues.

Elements of ${\mathcal W}(e,u)$ with respect to the same decomposition \eqref{H} look as
\eq{q1} q = a_{00}I  \oplus a_{01}I  \oplus a_{10}I  \oplus a_{11}I \oplus  \Phi(H), \en
where $\Phi=\begin{bmatrix} \phi_{00} & \phi_{01} \\ \phi_{10} & \phi_{11}\end{bmatrix}$, $a_{ij}\in\mathbb C$, and the functions $\phi_{ij}$ are Borel-measurable and essentially bounded on $\Delta$, in the sense of the spectral measure of $H$ (\cite{GiKu}, see also \cite{Spit94,BSpit10}). Consequently, $q\in{\mathcal W}(e,u)$ is an orthogonal projection if and only if $a_{ij}\in\{0,1\}$, the functions $\phi_{00}, \phi_{11}$ are real-valued, while $\phi_{01},\phi_{10}$ are complex conjugate, and
\eq{qpr} \phi_{00}-\phi_{00}^2= \phi_{11}-\phi_{11}^2=\abs{\phi_{01}}^2, \quad  (\phi_{00}+\phi_{11}-1)\phi_{01}=0. \en
On the other hand, direct computations immediately reveal that condition $quq=0$ is equivalent to
\eq{qu} \phi_{00}^2= \phi_{11}^2=\phi_{01}\phi_{10}, \quad  (\phi_{00}-\phi_{11})\phi_{01}=(\phi_{00}-\phi_{11})\phi_{10}=0. \en
Solving the system of equations \eqref{qpr}--\eqref{qu} yields
\[ \phi_{00}=\phi_{01}=\frac{1}{2}\chi, \ \phi_{01}=\frac{1}{2}\chi\omega, \ \phi_{10}=\frac{1}{2}\chi\overline{\omega} \]
with $\chi$ being a characteristic function of some subset of $\Delta$ and unimodular $\omega$.

So, elements of ${\mathcal Q}_u\cap{\mathcal W}(e,u)$ have the form
\eq{q} q = 0\oplus 0\oplus0\oplus0\oplus\frac{1}{2}\begin{bmatrix} \chi & \chi\omega \\ \chi\overline{\omega} & \chi\end{bmatrix}(H). \en

The rest of the reasoning depends on whether or not the subspaces ${\mathcal M}_{00}$, ${\mathcal M}_{11}$ are actually present in the decomposition \eqref{H}.

{\sl Case 1.} At least one of the subspaces ${\mathcal M}_{00}$, ${\mathcal M}_{11}$ is different from zero, that is, the range of $e$ contains some eigenvectors of $u$.

Since for any $q$ of the form \eqref{q} the restriction of $e-q$ on ${\mathcal M}_{00}\oplus{\mathcal M}_{11}$ is the identity, we then have $\norm{e-q}=1$. Consequently, $d=1$. Note that
in this case also $\norm{eue}=1$.

{\sl Case 2.} ${\mathcal M}_{00}={\mathcal M}_{11}=\{0\}$. Since both $e$ given by \eqref{tp} and $q$ given by \eqref{q} have zero restrictions onto ${\mathcal M}_{01}\oplus{\mathcal M}_{11}$,
we may without loss of generality suppose that in place of \eqref{H} simply ${\mathcal H}={\mathcal M}\oplus{\mathcal M}$, and respectively
\eq{eq1}  e = \begin{bmatrix}  H & \sqrt{H(I-H)} \\ \sqrt{H(I-H)} & I-H \end{bmatrix}, \quad q = \frac{1}{2}\begin{bmatrix} \chi & \chi\omega \\ \chi\overline{\omega} & \chi\end{bmatrix}(H). \en
So, $e-q=\Phi_{\chi,\omega}(H)$, where
\[ \Phi_{\chi,\omega}(t)=\begin{bmatrix}t-\frac{1}{2}\chi(t) & \sqrt{t(1-t)}-\frac{1}{2}\chi(t)\omega(t) \\  \sqrt{t(1-t)}-\frac{1}{2}\chi(t)\overline{\omega(t)} & 1-t-\frac{1}{2}\chi(t)\end{bmatrix}. \]
Consequently,
\[ \norm{e-q}=\es_{t\in\Delta}\lambda_{\chi,\omega}(t), \]
where $\lambda_{\chi,\omega}(t)$ is the positive eigenvalue of $\Phi_{\chi,\omega}(t)$, and {\em ess} is understood in the sense of the spectral measure of $H$.

If $\chi(t)=0$ for some $t\in\Delta$, then the respective $\lambda_{\chi,\omega}(t)$ equals one, guaranteeing $\norm{e-q}=1$. We should concentrate therefore on elements $q$ with $\chi(t)\equiv1$.
Then we have
\[ \Phi_{1,\omega}(t)=\begin{bmatrix}t-\frac{1}{2} & \sqrt{t(1-t)}-\frac{1}{2}\omega(t) \\  \sqrt{t(1-t)}-\frac{1}{2}\overline{\omega(t)} & \frac{1}{2}-t\end{bmatrix}, \]
and
\[ \lambda_{1,\omega}(t)=\sqrt{\frac{1}{2}-\sqrt{t(1-t)}\Re\omega(t)}.\]
Since $\omega$ is unimodular, to minimize $\lambda_{1,\omega}(t)$ for any given $t$ we should take $\omega(t)=1$. The respective element $q$ is simply \eq{q0} q_0=\frac{1}{2}\begin{bmatrix}I & I \\ I & I \end{bmatrix},\en
$\lambda_{1,1}(t)=\sqrt{\frac{1}{2}-\sqrt{t(1-t)}}$, and
\[ \norm{e-q_0}=\sqrt{\frac{1}{2}-\min_{t\in\Delta}\sqrt{t(1-t)}}=\sqrt{\frac{1}{2}\left(1-\sqrt{1-\max_{t\in\Delta}\abs{2t-1}^2}\right)}. \]
In order to justify \eqref{d}, it remains only to observe that
\eq{ne} \max_{t\in\Delta}\abs{2t-1}=\norm{eue}. \en But this is indeed the case, since
$eue=\Phi(H)$ with the matrix \[ \Phi(t)=(2t-1)\begin{bmatrix}t & \sqrt{t(1-1)} \\ \sqrt{t(1-1)} & 1-t\end{bmatrix},\] the eigenvalues of which are zero and $2t-1$.

\section{Additional comments}
\paragraph{1}
Recall \cite{VaSpit} that elements of $C^*$-algebra ${\mathcal C}(e,u)$ generated by $e$ and $u$ are those of the form \eqref{q1} for which the functions $\phi_{ij}$ are continuous on $\Delta$ and such that
$\phi_{01}(j)=\phi_{10}(j)=0$, $a_{ij}=\phi_{ii}(j)$ if $j\in\Delta$ ($i,j=0,1$). From \eqref{q0} we therefore conclude that the element $q_0\in{\mathcal W}(e,u)$ on which the distance from $e$ to ${\mathcal Q}_u$
is attained does not lie in ${\mathcal C}(e,u)$ if the spectrum of $H$ contains $0$ or $1$.

On the other hand, due to \eqref{ne} condition $\norm{eue}<1$ guarantees that $0,1\notin\Delta$, and thus the invertibility of the operator $H(I-H)$. Moreover, ${\mathcal M}_{00}={\mathcal M}_{11}=\{0\}$,
as was observed in Section~\ref{s:pr}. So, without loss of generality $e$ is given by the first formula in \eqref{eq1}, while $u=\dg [ I, -I ]$. From here:
\[  z:= \frac{1}{2}(e-ueu^*)=\begin{bmatrix}  0 & \sqrt{H(I-H)} \\ \sqrt{H(I-H)} & 0 \end{bmatrix}\in {\mathcal C}(e,ueu^*),\]
$z^2 = \dg [ H(I-H), H(I-H)]$ is positive definite and also lies in ${\mathcal C}(e,ueu^*)$, and therefore so does
$(z^2)^{-1/2} = \dg [\left((H(I-H)\right)^{-1/2}, \left(H(I-H)\right)^{-1/2}]$.
Along with $z$ and $(z^2)^{-1/2}$, the algebra ${\mathcal C}(e,ueu^*)$ contains their product
\[ z(z^2)^{-1/2}=\begin{bmatrix} 0 & I \\ I & 0 \end{bmatrix}.\]  We conclude from \eqref{q0} that $q_0=\frac{1}{2}(I+z(z^2)^{-1/2})\in {\mathcal C}(e,ueu^*)$.
\paragraph{2} The distances from $e$ to the sets ${\mathcal Q}_u$ and ${\mathcal Q}_u\cap{\mathcal W}(e,u)$ may not coincide. To illustrate, consider ${\mathcal H}={\mathbb C}^3$,
$u=\dg [1,1,-1]$ and $e=\dg [1,0,0]$. Then ${\mathcal Q}_u$ consists of zero and all matrices of the form
\[ q_{x,y}=\frac{1}{2}\begin{bmatrix} \abs{x}^2 & x\overline{y} & x \\ \overline{x}y &  \abs{y}^2& y \\ \overline {x} & \overline {y}  & 1\end{bmatrix}, \]
with the parameters $x,y\in\mathbb C$ satisfying $\abs{x}^2+\abs{y}^2=1$.
An easy computation shows that \[ \norm{e-q_{x,y}}=\sqrt{(1+\abs{y}^2)/2}. \]
So, the distance from $e$ to ${\mathcal Q}_u$ equals $1/\sqrt{2}$ and is attained on all the matrices $q_{\omega,0}$ with $\abs{\omega}=1$, that is, having the form
\[ \frac{1}{2}\begin{bmatrix} 1 & 0 & \omega \\ 0 &  0 & 0 \\  \overline {\omega}  &  0 &  1\end{bmatrix}. \]
On the other hand, the algebra generated by $e$ and $u$ consists simply of all $3$-by-$3$ diagonal matrices. The only diagonal matrix lying in ${\mathcal Q}_u$ is $0$,
and $d=1$ in full agreement with Theorem~\ref{th:main}.

\section{Acknowledgments} The author was supported in part by Faculty Research funding from the Division of Science and Mathematics, New York University Abu Dhabi.

\end{document}